\documentclass[a4paper,11pt]{amsart}

\usepackage{amsmath,amssymb,amscd}
\usepackage[arrow,matrix]{xy}
\usepackage{hyperref, pstricks}

\newtheorem{theorem}{Theorem}[section]

\newtheorem{lemma}{Lemma}[section]
\makeatletter
\renewcommand*{\c@lemma}{\c@theorem}
\makeatother

\newtheorem{proposition}{Proposition}[section]
\makeatletter
\renewcommand*{\c@proposition}{\c@theorem}
\makeatother

\newtheorem{corollary}{Corollary}[section]
\makeatletter
\renewcommand*{\c@corollary}{\c@theorem}
\makeatother

\makeatletter
\renewcommand*{\c@exercise}{\c@theorem}
\makeatother

\newtheorem{question}{Question}[section]
\makeatletter
\renewcommand*{\c@question}{\c@theorem}
\makeatother

\newtheorem{example}{Example}[section]
\makeatletter
\renewcommand*{\c@example}{\c@theorem}
\makeatother

\makeatletter
\renewcommand*{\c@conjecture}{\c@theorem}
\makeatother

\newtheorem{definition}{Definition}[section]
\makeatletter
\renewcommand*{\c@definition}{\c@theorem}
\makeatother

\makeatletter
\renewcommand*{\c@remark}{\c@theorem}
\makeatother

\newcommand{\Q}{\mathbb{Q}}

\newcommand{\mathscr}{}

\DeclareMathOperator{\rank}{rank}
\DeclareMathOperator{\vc}{vc}
\DeclareMathOperator{\Tor}{Tor}
\DeclareMathOperator{\Ker}{Ker}

\title[Toric cohomological rigidity of polytopes]{Toric cohomological rigidity of simple convex polytopes}

\author {Suyoung Choi}
\address{Department of Mathematical Sciences, KAIST, 335 Gwahangno, Yuseong-gu, Daejeon 305-701, Republic of
Korea,
\newline{\it and}\newline
(current address)
   Department of Mathematics, Osaka City University,
   Sugimoto, Sumiyoshi-ku, Osaka 558-8585, Japan}

\email{choisy@kaist.ac.kr\ \ \ choi@sci.osaka-cu.ac.jp}
\author{Taras Panov}
\address{Department of Geometry and Topology, Faculty od Mathematics
and Mechanics, Moscow State University, Leninskiye Gory, Moscow
119991, Russia,\newline {\it and}\newline Institute for
Theoretical and Experimental Physics, Moscow 117259, Russia}
\email{tpanov@mech.math.msu.su}
\author{Dong Youp Suh}
\address{Department of Mathematical Sciences, KAIST, 335 Gwahangno, Yuseong-gu, Daejeon 305-701, Republic of Korea}
\email{dysuh@math.kaist.ac.kr}

\subjclass{55Nxx, 52Bxx}

\thanks{The first author was partially supported by Brain
Korea 21 project, KAIST; the second author was partially supported
by the Russian Foundation for Basic Research (grants
no.~09-01-00142, 08-01-91855-KO), the State Programme for the
Support of Leading Scientific Schools (grant 1824.2008.1), and
P.~Deligne's 2004 Balzan prize in mathematics; the third author
was partially supported by Basic Science Research Program through
the National Research Foundation of Korea (NRF) funded by the
Ministry of Education, Science and Technology (2009-0063179 )}

\begin{document}
\begin{abstract}
A simple convex polytope $P$ is \emph{cohomologically rigid} if
its combinatorial structure is determined by the cohomology ring
of a quasitoric manifold over $P$. Not every $P$ has this
property, but some important polytopes such as simplices or cubes
are known to be cohomologically rigid. In this article we
investigate the cohomological rigidity of polytopes and establish
it for several new classes of polytopes including products of
simplices. Cohomological rigidity of $P$ is related to the
\emph{bigraded Betti numbers} of its \emph{Stanley--Reisner ring},
another important invariant coming from combinatorial commutative
algebra.
\end{abstract}

\maketitle

\section{Introduction}

Quasitoric manifolds were defined by Davis and Januszkiewicz in
\cite{DJ} as a topological analogue of nonsingular toric
varieties. Namely, a \emph{quasitoric manifold} over a simple
convex polytope $P$ is a closed $2n$-dimensional manifold $M$ with
a locally standard action of an $n$-torus $G=(S^1)^n$ (that is,
the action locally looks like a faithful real $2n$-dimensional
representation of $G$) and a surjective map $\pi\colon M\to P$
whose fibers are the $G$-orbits. The combinatorial structure of
$P$ is completely determined by the equivariant cohomology ring
$H^\ast_G(M)$ because  $H^\ast_G(M)$ is isomorphic to the
Stanley-Reisner ring (or the face ring) of $P$. On the other hand
the $2i$-th
Betti number of $M$ is equal to the $i$-th component of the
$h$-vector of $P$. Therefore the usual cohomolgy $H^\ast(M)$
contains some combinatorial information of~$P$.

In general the cohomology ring of a quasitoric manifold does not
contain sufficient information to determine the combinatorial
structure of the base polytope $P$, as in Example~4.3 of
\cite{MS}, which we will discuss briefly for reader's convenience.
To do this let us fix some notation. For an $n$-dimensional simple
convex polytope $P$ and a vertex $v$ of it, let $\vc(P, v)$ denote
the connected sum of $P$ with the $n$-simplex $\Delta^n$ at the
vertex $v$. Hence $\vc(P, v)$ is the simple convex polytope
obtained from $P$ by cutting a small $n$-simplex neighborhood of
the vertex $v$. We call $\vc(P, v)$ the {\em vertex cut} of $P$ at
$v$. When the combinatorial structure of $\vc(P, v)$ does not
depend on the vertex $v$, we simply denote it  by $\vc(P)$. For
example when $P$ is a product of simplices, the vertex cut
$\vc(P,v)$ does not depend on the choice of a  vertex~$v$.

The following example explains a phenomenon leading to our main
definition.

\begin{example}\label{vcex}
We
consider $M=\mathbb CP^2\times \mathbb CP^1$ with the standard
$(S^1)^3$-action. It is a quasitoric manifold over the triangular
prism $P=\Delta^2\times \Delta^1$. The equivariant blow up $M'$ of
$M$ at a fixed point $x$ is a quasitoric manifold over
$P'=\vc(P)$, which does not depend on the choice of a fixed
point~$x$. Now if we blow up $M'$ equivariantly at a fixed point
$y$ in $M'$, then the resulting manifold $M''$ is a quasitoric
manifold over $P''=\vc(P',v)$. The manifold $M''$ is no longer
independent of a fixed point $y$; in fact there are three
equivariantly different manifolds corresponding to three
combinatorially different vertex cuts $\vc(P',v)$ (these
correspond to the first three simplicial complexes in the second
line in p.~192 of~\cite{Od}).

On the other hand, the cohomology ring of $M''$ does not depend on
the choice of a fixed point $y$, because $M''$ is homeomorphic to
the connected sum of $\mathbb CP^2\times \mathbb CP^1$ with two
copies of $\mathbb CP^3$.
We therefore are in the situation when the cohomology ring of a
quasitoric manifold does not determine the combinatorial structure
of the base polytope.
\end{example}

Nevertheless, in many cases the combinatorial type of $P$ is
determined by $H^\ast(M)$. We therefore naturally come to the
following definition, firstly introduced in~\cite{MS}.

\begin{definition}\label{definiton:rigidity}
A simple polytope $P$ is \emph{cohomologically rigid} if  there exists a
quasitoric manifold $M$ over $P$, and whenever there exists a
quasitoric manifold $N$ over another polytope $Q$
with a graded ring isomorphism $H^\ast(M)\cong H^\ast(N)$ there is
a combinatorial equivalence $P\approx Q$. We shall refer to
such $P$ simply as \emph{rigid} throughout the paper.
\end{definition}

We shall extend this definition to arbitrary Cohen-Macaulay
complexes in Definition~\ref{definition:cohen-macauley complex}.
In~\cite{MS} the rigidity property is expressed in terms of toric
manifolds, but here we modify the original definition to make use
of a wider class of quasitoric manifolds.
The interval $I$ is trivially rigid. More
generally,
it is shown in~\cite{MP} that any cube $I^n$ is rigid. In
Section~2 we give more classes of rigid polytopes, as described in
the following results.

\setcounter{section}2 \setcounter{theorem}1

\begin{theorem}\label{theorem:only one polytope}
Let $P$ be a simple polytope supporting a quasitoric manifold.
If there is no other simple polytope
with the same numbers of $i$-faces as those of $P$ for all $i$,
then $P$ is rigid.
\end{theorem}

\begin{corollary}\label{corollary: polygon is rigid}
Every polygon, i.e. 2-dimensional convex polytope, is rigid.
\end{corollary}

A  simple convex polytope is called {\em triangle-free} if it  has no triangular $2$-face.
The following result, proved in Section~4, establishes the
rigidity for triangle-free polytope with few facets.

\setcounter{section}4 \setcounter{theorem}2

\begin{theorem}\label{theorem: trianlgefree polytope}
Every triangle-free $n$-dimensional simple convex polytope with less than $2n+3$ facets
is rigid.
\end{theorem}

Since a cube $I^n$ has $2n$ facets, \autoref{theorem: trianlgefree
polytope} gives a different proof of the rigidity for cubes from
that in ~\cite{MP}.
In Section~5 this result is generalized as follows

\setcounter{section}5 \setcounter{theorem}2

\begin{theorem}\label{theorem: product of simplices is rigid}
A finite product of simplices is rigid.
\end{theorem}

From the argument in Example~\ref{vcex} one can see immediately
that if the vertex cut of a polytope $P$ depends on a choice of
vertex, then all the vertex cuts of $P$ are not rigid. So it is
natural to ask whether $\vc(P)$ is rigid if the vertex cut of $P$
is independent of a choice of vertex.
In section~6 we confirm
this when $P$ is a product of simplices:

\setcounter{section}6 \setcounter{theorem}3

\begin{theorem}\label{theorem: cv of product of simplices}
If $P$ is a finite product of simplices, then $\vc(P)$ is  rigid.
\end{theorem}

\setcounter{section}1 \setcounter{theorem}2

We can apply the above results to determine rigidity of  $3$-dimensional simple convex polytopes
with facet numbers up to nine.  This result is given in \autoref{section:3-polytope}. We also prove
that dodecahedron is rigid in \autoref{theorem:dodecahedron}.

The rigidity property for simple polytopes is closely related to
the following interesting question on quasitoric manifolds.

\begin{question}\label{q1}
Suppose $M$ and $N$ are two quasitoric manifolds such that
$H^\ast(M)\cong H^\ast(N)$ as graded rings. Are $M$ and $N$
homeomorphic?
\end{question}

We can also consider the following slightly weaker question, which
can be considered as an intermediate step to answering
Question~\ref{q1}.

\begin{question}\label{q2}
Suppose $M$ and $N$ are two quasitoric manifolds over the same
simple convex polytope $P$ such that $H^\ast(M)\cong H^\ast(N)$ as
graded rings. Are $M$ and $N$ homeomorphic?
\end{question}

Question~\ref{q1} for quasitoric manifolds whose cohomology rings
are isomorphic to those of a product of copies of $\mathbb CP^1$ is
considered in~\cite{MP}, and it is shown there that these
manifolds are actually homeomorphic to a product of copies of $\mathbb
CP^1$. This is done in two steps; firstly the result is proved
under additional assumption that the quotient polytope is a cube
$I^n$, and then the rigidity of $I^n$ is established,
see~\cite{MP}.

In \cite{CMS} it is proved that if $M$ is a quasitoric manifold  over a product of simplices
$\prod_{i=1}^t\Delta^{n_i}$ such that $H^\ast(M)\cong H^\ast(\prod_{i=1}^t\mathbb CP^{n_i})$,
then $M$ is homeomorphic to $\prod_{i=1}^t\mathbb CP^{n_i}$. Since a product of simplices is rigid by
\autoref{theorem: product of simplices is rigid}, we have the following theorem.

\begin{theorem}
Suppose $M$ is a quasitoric manifold such that $H^\ast(M)\cong H^\ast(\prod_{i=1}^t\mathbb CP^{n_i})$ as graded rings,
then $M$ is homeomorphic to $\prod_{i=1}^t\mathbb CP^{n_i}$.
\end{theorem}

The main technical ingredient for the proofs of the results in
this paper is the following proposition. For a polytope $P$ let
$\beta^{-i, 2j}(P)$ be the bigraded Betti numbers of the
Stanley-Reisner ring $\mathbb Q(P)$ of $P$, see \autoref{section:
bigraded Betti number} or \cite{BP} for details.

\setcounter{section}3 \setcounter{theorem}7

\begin{proposition}\label{proposition: bigraded Betti number}
Let $M$ (reps. $N$) be a quasitoric manifold over $P$ (resp. $Q$).
If $H^\ast(M : \mathbb Q)\cong H^\ast(N:\mathbb Q)$ as graded rings, then
$\beta^{-i, 2j}(P)=\beta^{-i, 2j}(Q)$ for all $i$ and $j$.
\end{proposition}

\setcounter{section}1

\section{Rigidity and $f$-vectors}
For a convex $n$-dimensional polytope $P$ let $f_i$ denote the
number of codimension  $i+1$ faces of $P$, and let $f(P)=(f_0,
\ldots, f_{n-1})$ denote the $f$-vector of $P$. Note that if $P$
and $Q$ are two $2$-dimensional polytopes, then $f(P)=f(Q)$
implies $P\approx Q$. Recall that $h$-vector $h(P)=(h_0, \ldots,
h_n)$ of $P$ is defined by
$$\sum_{i=0}^n h_i t^{n-i}=\sum_{j=0}^n f_{j-1}(t-1)^{n-j}.$$
The following theorem  proved in \cite{DJ} shows that $f$-vector  of the base polytope
$P$ is determined by
the cohomology ring of the quasitoric manifold $M$ over $P$.

\begin{theorem}[\cite{DJ}]
\label{theorem: h_i = 2i-th Betti number} For a quasitoric
manifold $M$ over $P$ the $2i$-th Betti number $b_{2i}(M)$ of $M$
is equal to the $i$-th component $h_i$ of the $h$-vector of $P$.
\end{theorem}


\begin{theorem}
Let $P$ be a
simple polytope supporting a quasitoric manifold. If there is no
other simple polytope with the same numbers of $i$-faces as those
of $P$ for all $i$, then $P$ is rigid.
\end{theorem}

\begin{proof}
Now let $P$ be a polytope and $M$ a given quasitoric manifold over
$P$. Suppose $N$ is another quasitoric manifold over $Q$ such that
$H^\ast(M)\cong H^\ast(N)$ as graded rings. Then the cohomology
isomorphism implies $b_{2i}(M)=b_{2i}(N)$ for all $i$. Hence
$h(P)=h(Q)$ by \autoref{theorem: h_i = 2i-th Betti number}, which
implies $f(P)=f(Q)$. Since there is no other simple convex
polytope with the same face numbers of $P$, $P \approx Q$ and
hence \autoref{theorem:only one polytope} is proved.
\end{proof}

\begin{corollary}
Every polygon, i.e. 2-dimensional convex polytope, is rigid.
\end{corollary}

\begin{proof}
\autoref{corollary: polygon is rigid} follows immediately from
\autoref{theorem:only one polytope}.
\end{proof}

\section{Bigraded Betti numbers of polytopes}\label{section: bigraded Betti number}

Let $A=\mathbb Q[x_1,\ldots, x_m]$ be the polynomial graded ring in
$x_1,\ldots, x_m$ over the rationals with $\deg x_i=2$ for all $i$.
A {\em free resolution} $[R,d]$ of a finitely generated $A$-module $M$
is an exact sequence
\begin{equation}\label{equation: free resolution}
\xymatrix{ 0 \ar[r]& R^{-n} \ar[r]^-{d}& R^{-n+1} \ar[r]^-{d}&
\cdots\ar[r]^-{d}& R^0\ar[r]^-{d}& M \ar[r] & 0},
\end{equation}
where $R^{-i}$ are finitely generated free graded $A$-modules and
$d$ are degree preserving homomorphisms. If we take $R^{-i}$ to be
the module generated by the minimal basis of $\Ker(d:R^{-i+1}\to
R^{-i+2})$, then we get a {\em minimal resolution} of $M$. This also
shows the existence of a resolution.

Dropping the last term $M$ in the sequence~\eqref{equation: free
resolution} and tensoring it over $A$ with another finitely
generated $A$-module $N$, we obtain the following sequence:

\begin{equation}\label{equation: tensored resolution}
\xymatrix{ 0 \ar[r]& {R^{-n}\otimes_AN} \ar[r]^-{d\otimes 1}& {R^{-n+1}\otimes_AN} \ar[r]^-{d\otimes 1}&
\cdots\ar[r]^-{d\otimes 1}& R^0\otimes_AN\ar[r] & 0}
\end{equation}
This sequence is not necessarily exact, and its cohomology is
known as the \emph{Tor-modules}:
$$\Tor_A^{-i}(M,N):=H^{i}(R^{-\ast}\otimes_AN).$$
Since everything is graded, we actually have the grading
$$\Tor_A^{-i}(M,N)=\bigoplus_j\Tor_A^{-i, j}(M,N).$$
The following proposition is well-known, and we refer the reader to \cite{BP}
for details.

\begin{proposition}\label{proposition: properties of Tor}
The above defined Tor-modules satisfy the following properties.
\begin{enumerate}
\item $\Tor_A(M,N)$ is independent of the choice of a resolution $[R,d]$ of $M$.
\item $\Tor_A(M,N)$ is functorial in all three arguments, i.e., in $A$, in $M$, and in $N$.
\item $\Tor_A^0(M,N)=M\otimes_A N$.

\item $\Tor_A^{-i}(M,N)=\Tor_A^{-i}(N,M)$.
\end{enumerate}
\end{proposition}
\begin{flushright}
    $\Box$
\end{flushright}

We regard $\mathbb Q$ as an $A$-module via the ring map $A\to
\mathbb Q$ sending each $x_i$ to~$0$.
Set $N=\mathbb Q$ and consider $\Tor_A(M,\mathbb Q)$.
\begin{definition}
The {\em bigraded Betti numbers} of $M$ are defined by
$$\beta^{-i}(M)=\dim_{\mathbb Q}\Tor_A^{-i}(M,\mathbb Q), $$
$$\beta^{-i,j}(M)=\dim_{\mathbb Q}\Tor_A^{-i,j}(M,\mathbb Q).$$
\end{definition}

When $[R,d]$ is a minimal resolution of $M$, then the map
$$d\otimes 1:R^{-i}\otimes_A\mathbb Q\to R^{-i+1}\otimes_A\mathbb Q$$
are the zero maps for all $i$.
Hence $\beta^{-i, j}=\rank_\Q R^{-i,j}$.

We now consider the case when $M$ is the Stanley-Reisner ring $\mathbb Q(P)$ of a
simple convex polytope $P$, which is
$$\mathbb Q(P)=\mathbb Q[x_1,\ldots, x_m]/I_P$$
where $x_i$ are indeterminates corresponding to the facets $F_i$
of $P$, $m$ is the number of facets, and $I_P$ is the homogeneous
ideal generated by the monomials $x_{i_1}\cdots x_{i_\ell}$
whenever $F_{i_1}\cap\cdots\cap F_{i_\ell}=\emptyset$. This $I_P$
is called the {\em Stanley-Reisner ideal} of $P$. Then $\mathbb
Q(P)$ is a graded $A$-module with $\deg x_i=2$ for all
$i=1,\ldots, m$. The {\em bigraded Betti numbers} of $P$
are defined to be $\beta^{-i,2j}(P)=\beta^{-i,2j}(\mathbb Q(P))$.
Since $\deg x_i=2$ we only have even index $2j$.

>From the previous observation that $\beta^{-i, j}=\rank_\mathbb Q
R^{-i,j}$ for a minimal resolution $[R,d]$, we can see easily that
$\beta^{-1, 2j}$ is equal to the number of degree $2j$ monomial
elements in a minimal basis of the ideal $I_P$. For example, if
$P=I^n$ then $x_ix_{n+i}$ for $i=1,\ldots n$ form a minimal basis
for the Stanley-Reisner ideal $I_P$ of $P$
(here we assume that $x_i$ and $x_{n+i}$ are the generators
corresponding to the opposite facets $F_i$ and $F_{n+i}$ of
$I^n$). Hence
$$\beta^{-1, 2j}(I^n)=\begin{cases} n,  & j=2;\\
                                    0,  & \textrm{otherwise}.\end{cases}$$

The following theorem of Hochster gives a nice formula
for bigraded Betti numbers.

\begin{theorem}[\cite{Ho}]
\label{theorem: Hochster} Let $P$ be a simple convex polytope with
facets $F_1$,$\ldots$,$F_m$. For a subset
$\sigma\subset\{1,\ldots,m\}$ let
$P_\sigma=\cup_{i\in\sigma}F_i\subset P$. Then we have
$$\beta^{-i, 2j}(P)=\sum_{|\sigma|=j}\dim\widetilde H^{j-i-1}(P_\sigma).$$
Here  $\dim\widetilde H^{-1}(\emptyset)=1$ by convention.
\end{theorem}

Bigraded Betti numbers also satisfy the following relations, see
\cite{BP} for details.

\begin{proposition}\label{propositon: relations of Betti numbers}
Let $P$ be an $n$-dimensional simple convex polytope with $m$ facets, i.e., $f_0(P)=m$.
Then
\begin{enumerate}
\item $\beta^{0,0}(P)=\beta^{-(m-n), 2m} (P)=1$,
\item (Poincar\'{e} duality) $\beta^{-i, 2j}(P)=\beta^{-(m-n)+i, 2(m-j)}(P)$, and
\item $\beta^{-i,2j}(P_1\times P_2)=\sum_{i'+i''=i, j'+j''=j}\beta^{-i',2j'}(P_1)
\beta^{-i'', 2j''}(P_2)$.
\end{enumerate}
\end{proposition}
\begin{flushright}
    $\Box$
\end{flushright}

\begin{definition}
A sequence $\lambda_1,\ldots,\lambda_p$ of homogeneous elements in
$\mathbb Q(P)$ is a {\em regular sequence} if it is algebraically
independent and $\mathbb Q(P)$ is a free module over $\mathbb
Q[\lambda_1,\ldots, \lambda_p]$.
\end{definition}

Let $J$ be an ideal of $\mathbb Q(P)$ generated by a regular sequence $\lambda_1, \ldots, \lambda_p$.
Let $\pi:A\to \mathbb Q(P)$ be the projection.
Choose
homogeneous
$t_i\in A$ such that $\pi(t_i)=\lambda_i$. Let $J$
also denote the ideal of $A$ generated by $t_1, \ldots, t_p$.

\begin{lemma}[{\cite[Lemma 3.35]{BP}}]
\label{lemma: Tor with regular sequence} Let $J$ be an ideal
generated by a regular sequence of $\mathbb Q(P)$. Then we have
the following algebra isomorphism.
$$\Tor_A^{\ast, \ast}(\mathbb Q(P), \mathbb Q)\cong \Tor_{A/J}^{\ast, \ast}(\mathbb Q(P)/J, \mathbb Q)$$
\end{lemma}
\begin{flushright}
    $\Box$
\end{flushright}

\begin{lemma}\label{lemma: Tor(P)=Tor(p')}
Let $P$ and $P'$ be two $n$-dimensional simple convex polytopes.
Let $J=(\lambda_1, \ldots, \lambda_n)$ (resp. $J'=(\lambda'_1,
\ldots, \lambda'_n)$) be an ideal of $\mathbb Q(P)$ (resp.
$\mathbb Q(P')$) generated by a regular sequence of  degree $2$
elements $\lambda_i$ (resp. $\lambda'_i$).
If there is a graded ring isomorphism $h\colon\mathbb
Q(P)/J\stackrel{\cong}\longrightarrow \mathbb Q(P')/J'$, then
$f_0(P)=f_0(P')$ and
$$\Tor_A^{\ast, \ast}\bigl(\mathbb Q(P), \mathbb Q\bigr)=
\Tor_A^{\ast, \ast}\bigl(\mathbb Q(P'),\mathbb Q\bigr).$$
\end{lemma}

\begin{proof}
Note that the Stanley-Reisner ring $\mathbb Q(P)$ is generated by
$f_0(P)$
elements of degree two. Since $J$ and $J'$ are generated by degree
two elements, the equality $f_0(P)=f_0(P')$ follows immediately
from the isomorphism of degree two subgroups induced from $\mathbb
Q(P)/J\cong \mathbb Q(P')/J'$. Thus we may assume that $J$ and
$J'$ are both ideals of $A=\mathbb Q[x_1,\ldots,x_m]$ and $h$ is
an $A$-algebra isomorphism. By \autoref{lemma: Tor with regular
sequence} we have $\Tor_A(\mathbb Q(P), \mathbb
Q)=\Tor_{A/J}(\mathbb Q(P)/J,\mathbb Q)$, and a similar equality
holds for $P'$.

Now we claim that there is an $A$-algebra isomorphism $\bar
h\colon A/J\to A/J'$ closing the commutative diagram
\begin{equation}\label{torfu}
\begin{CD}
  A/J @>{\bar h}>> A/J'\\
  @VVV @VVV\\
  \mathbb Q(P)/J @>h>> \mathbb Q(P')/J'.
\end{CD}
\end{equation}
Note that both $A/J$ and $A/J'$ are isomorphic to $\mathbb
Q[x_1,\ldots, x_{m-n}]$ where $m=f_0(P)=f_0(P')$.  Also note that
the projection maps $A/J\to \mathbb Q(P)/J$ and $A/J'\to \mathbb
Q(P')/J'$ induce isomorphisms $(A/J)_2\to (\mathbb Q(P)/J)_2$ and
$(A/J')_2\to (\mathbb Q(P')/J')_2$ on degree $2$ subgroups.
Therefore we have an isomorphism $(A/J)_2\to (\mathbb Q(P)/J)_2\to
(\mathbb Q(P')/J')_2\to (A/J')_2$. Since $A/J$ and $A/J'$ are
generated in degree~2, we obtain the isomorphism $\bar h\colon
A/J\to A/J'$ as necessary.

Finally, the required isomorphism
\[
  \Tor_{A/J}^{\ast, \ast}\bigl(\mathbb Q(P)/J, \mathbb Q\bigr) \cong
  \Tor_{A/J'}^{\ast, \ast}\bigl(\mathbb Q(P')/J', \mathbb Q\bigr)
\]
follows from~\eqref{torfu} and the functoriality of $\Tor$ in
\autoref{proposition: properties of Tor}(2).
\end{proof}

We are now ready to prove
the invariance of the bigraded Betti numbers.

\begin{proposition}
Let
$M$ (reps. $N$) be a quasitoric manifold over $P$ (resp. $Q$). If
$H^\ast(M : \mathbb Q)\cong H^\ast(N:\mathbb Q)$ as graded rings,
then $\beta^{-i, 2j}(P)=\beta^{-i, 2j}(Q)$ for all $i$ and $j$.
\end{proposition}

\begin{proof}
Recall that if $M$ is a quasitoric manifold over a simple convex polytope $P$, then
$H^\ast(M\colon \mathbb Q)\cong \mathbb Q[x_1, \dots, x_m]/K$ where $K=I_P +J$ and $I_P$ is the
rational Stanley-Reisner ideal of $P$, and $J$ is an ideal generated by some linear combinations
$\lambda_{i1}x_1+\dots+ \lambda_{im}x_m \in \mathbb Q[x_1,\dots x_m]$ for $i=1,\dots, n$ which
project to a regular sequence $\theta_1, \dots, \theta_n$ in $\mathbb Q[x_1,\dots, x_m]/I_P$, see
\cite{DJ}.
Here $m$ is the number of facets in $P$.
Therefore we have
the isomorphism
$$\mathbb Q(P)/J\cong H^\ast(M:\mathbb Q)\cong H^\ast(N:\mathbb Q)\cong \mathbb Q(P')/J'.$$
Hence the proposition follows from \autoref{lemma: Tor(P)=Tor(p')}.
\end{proof}

Since $\beta^{-i,j}=\rank_\mathbb QR^{-i,j}$ for a minimal
resolution $[R,d]$, and since
the Betti numbers are independent of the choice of a resolution,
it is convenient to calculate $\beta^{-i,j}$ using a particular
minimal resolution. For this purpose we will consider the minimal
resolution of $\mathbb Q(P)$  corresponding to the {\em canonical
minimal basis} of the rational Stanley-Reisner ideal $I_P$, which
we define below. The following procedure is explained in Example
3.2 in\cite{BP};
we also reproduce it here for the reader's convenience.

In general, for a finitely generated graded $A$ module $M$ the canonical minimal basis can
be chosen as follows. Take the lowest degree, say $d_1$, elements in $M$ which
form a $\mathbb Q$-vector subspace of $M$, and choose its basis $\mathcal B_{d_1}$.
Then span an $A$-submodule $M_1$ of $M$ spanned by $\mathcal B_{d_1}$.
Then take the lowest degree, say $d_2$, elements in $M\setminus M_1$ which form
a $\mathbb Q$-subspace of $M$, and choose its basis $\mathcal B_{d_2}$.
Span an $A$-submodule $M_2$ of $M$ with $\mathcal B_{d_1}\cup\mathcal B_{d_2}$.
Continue this process. Since $M$ is finitely generated, this process must stop at some
$p$-th step, and we get $\mathbb Q$-subspace $\mathcal B_{d_p}$ for $M\setminus M_{p-1}$.
Then $M$ is generated by $\mathcal B=\cup_{i=1}^p\mathcal B_{d_i}$ as an $A$-module.
The generator set $\mathcal B$ constructed in this way has
the minimal possible number of elements, and we call it the
canonical minimal basis of $M$.

In particular,
for $M=I_P$ the canonical basis is $\mathcal B=\cup_{\ell\ge
1}\mathcal B_{2\ell}$ where $\mathcal B_{2\ell}$ are inductively
defined as follows. $\mathcal B_{2}$ consists of all monomials
$x_ix_j$ such that $F_{i}\cap F_{j}= \emptyset$ where $F_k$ is the
facet of $P$ corresponding to $x_k$. Assume $\mathcal B_{2k}$ is
defined for $ k<\ell$. Then $\mathcal B_{2\ell}$ consists of the
monomials $x_{i_1}\cdots x_{i_{\ell}}$ that are not divisible by
the elements in $\cup_{i=1}^{\ell-1}\mathcal B_{2i}$ such that
$\cap_{k=1}^\ell F_{i_k}=\emptyset$

For finitely generated $A$ module $N$, there is the following way
of constructing minimal resolution of $N$. Take a minimal basis
$\mathcal B_N$,
and define $R^0$ to be a free $A$ module
generated by the elements of
$\mathcal B_N$. There is an obvious epimorphism $R^0\to N$.  Take
a minimal basis
for $\ker(R^0\to M)$, and define $R^{-1}$ to be a free $A$ module
with these generators, and so on.


\begin{example}
1. If $P=I^n$, the $n$-dimensional cube, then
$${\mathcal B}(I_P)
=\{x_i x_{n+i}\mid i=1,\ldots, m\}.$$

2. If $P=\prod_{i=1}^t\Delta^{n_i}$, a product of simplices, then
$${\mathcal B}(I_P)=\{x_{i,0}\cdots x_{i,
n_i}\mid i=1,\ldots, t\}.$$
\end{example}

We close this section by giving an algebraic version of rigidity.
Recall that the rational Stanley-Reisner ring $\mathbb Q(K)$ of a
simplicial complex $K$ with $m$ vertices $v_1,\ldots, v_m$ is the
quotient ring $\mathbb Q[x_1,\ldots, x_m]/I_K$ where $I_K$ is the
ideal generated by the monomials $x_{i_1}\cdots x_{i_\ell}$ where
the corresponding vertices $v_{i_1}, \ldots, v_{i_\ell}$ do not
form a simplex on $K$. Then the rational Stanley-Reisner ring
$\mathbb Q(P)$ of a simple convex polytope $P$ is actually the
rational Stanley-Reisner ring of the dual simplicial complex of
$\partial P$, i.e., $\mathbb Q(P)=\mathbb Q((\partial P)^\ast)$.
Since $P$ is simple $(\partial P)^\ast$ is a simplicial complex.

The above constructed minimal basis $\mathcal B$ of $I_P$
coincides with the canonical minimal basis of the ideal $I_K$
(see~\cite[\S3.4]{BP}) consisting of monomials corresponding to
all \emph{missing faces} of the simplicial complex $K$ dual to the
boundary of $P$ (a missing face of a simplicial complex is its
subset of vertices which does not span a simplex, but every whose
proper subset does span a simplex).

A simplicial complex of dimension $n-1$ is called {\em
Cohen-Macaulay} if there exists a length $n$ regular sequence in
$\mathbb Q(K)$. For any $n$-dimensional simple convex polytope
$P$, its dual $(\partial P)^\ast$ is known to be  Cohen-Macaulay.
Therefore the definition of rigidity of a simple polytope
can be generalized to that of a Cohen-Macaulay complex as follows:

\begin{definition}\label{definition:cohen-macauley complex}
An $(n-1)$-dimensional Cohen-Macaulay complex $K$ is {\em rigid}
if for any $(n-1)$-dimensional Cohen-Macaulay complex $K'$ and for
ideals $J\subset \mathbb Q(K)$ and $J'\subset \mathbb Q(K')$
generated by degree $2$ regular sequences of length $n$, $\mathbb
Q(K)/J\cong\mathbb Q(K')/J'$ implies $\mathbb Q(K)\cong\mathbb
Q(K')$.
\end{definition}

\section{Rigidity of triangle-free simple polytopes}
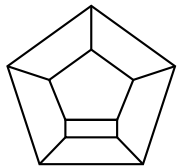
\begin{figure}
\psset{unit=11pt}
\begin{pspicture}(0,0)(6,6)
\pspolygon(3,6)(0.15,3.93)(1.24,0.57)(4.76,0.57)(5.85,3.93)
\pspolygon(3,4.5)(1.57,3.46)(2.12,2.09)(3.88,2.09)(4.43,3.46)
\pspolygon(2.12,2.09)(3.88,2.09)(3.88,1.49)(2.12,1.49)
\psline(3,6)(3,4.5)
\psline(0.15,3.93)(1.57,3.46)
\psline(1.24,0.57)(2.12,1.49)
\psline(4.76,0.57)(3.88,1.49)
\psline(5.85,3.93)(4.43,3.46)
\end{pspicture}
  \caption{Schlegel diagram of $Q$}\label{figure:4^45^4}
\end{figure}

It is shown in \cite{BB90} that if $P$ is a triangle-free convex
$n$-polytope then $f_i(P)\ge f_i(I^n)$ for all $i=0,\ldots, n-1$.
Therefore the number of facets of $P$ satisfies $f_{0}(P)\ge 2n$.
Furthermore it is shown in \cite{BB92} that if $P$ is simple and
\begin{enumerate}
\item if $f_{0}(P)=2n$, then $P\approx I^n$,
\item if $f_{0}(P)=2n+1$, then $P\approx P_5\times I^{n-2}$
where $P_5$ is a pentagon, and
\item if $f_0=2n+2$, then $P  \approx P_6\times I^{n-2}$, $Q\times I^{n-3}$, or
$P_5\times P_5\times I^{n-4}$ where $P_6$ is the hexagon and $Q$ is $3$-dimensional simple convex polytope
obtained from pentagonal prism by cutting out one of the edges forming a pentagonal facet, see \autoref{figure:4^45^4}.
\end{enumerate}

\begin{lemma}\label{lemma:minimal basis and triangle-free}
Let $P$ be an $n$-dimensional simple
polytope.
If $\beta^{-1, 2j}(P)=0$ for all $j\ge 3$, then $P$ is triangle-free.
\end{lemma}

\begin{proof}
Suppose otherwise. Namely, suppose there exists a triangular
$2$-face $T$ of $P$. Each edge $e_i$ of $T$ for $i=1,2,3$ is
an intersection of $n-1$  facets of $P$.  Thus there exists a
unique facet, say $F_i$ which contains the edge $e_i$ but not the
triangle $T$ for $i=1,2,3$.
Since $P$ is simple, $T$ is the intersection of exactly $n-2$
facets, and we may assume that $T=\cap_{i=4}^{n+1}F_i$.
Then $\cap_{i=1}^{n+1}F_i=\emptyset$ because $P$ is simple. This
means that the monomial $\prod_{i=1}^{n+1}x_i$ is contained in
$I_P$, where $x_i$ is the degree two generating element of $I_P$
corresponding to the facet $F_i$ for $i=1,\ldots n+1$. Therefore
there exists a minimal basis element $x_{n_1}\cdots x_{n_k}$ of
$I_P$ that divides $\prod_{i=1}^{n+1}x_i$. Now consider the set
$\mathcal S=\{F_{n_1}, \ldots,
F_{n_k}\}$
of facets corresponding to $x_{n_i}$ for $i=1,\ldots, k$. Then the
intersection the elements of any proper subset of $\mathcal S$ is
nonempty, but the intersection of the elements of $\mathcal S$ is
empty. Note that $\cap_{i=1,i\ne j}^{n+1}F_i=v_j$ for $j=1,2,3$
where $v_j$ is the opposite vertex of $T$ to the edge $e_j$.
Therefore $\mathcal S$ must contain the facets $F_1$, $F_2$, and
$F_3$. Therefore the minimal basis element $x_{n_1}\cdots x_{n_k}$
should divide $x_1x_2x_3$, hence $x_{n_1}\cdots x_{n_k}$ is of
degree greater than or equal to $6$ which  contradicts to the
hypothesis $\beta^{-1, 2j}(P)=0$ for all $j \ge 3$.
%
%
\end{proof}
Note that the condition $\beta^{-1, 2j}(P)=0$ for all $j\ge 3$ means that
the Stanley-Reisner ideal $I_P$ of $P$ is generated by quadratic monomials
of the form $x_ix_j$,
and this is equivalent to saying that the simplicial complex $K=(\partial P)^\ast$
is flag.

If the number of facets of $P$ is less than or equal to $2n+2$, then the converse of
\autoref{lemma:minimal basis and triangle-free} is true. Namely, we have
\begin{lemma}\label{lemma: necessatity for triangle-free}
If $P$ is a triangle-free $n$-dimensional simple convex polytope with $f_0(P)\le 2n+2$, then
$\beta^{-1, 2j}(P)=0$ for all $j\ge 3$.
\end{lemma}

\begin{proof}
Since $f_0(P)\le 2n+2$, we know that $P\approx I^n$, $P_5\times I^{n-2}$, $P_6\times I^{n-2}$,
$Q\times I^{n-3}$, or $P_5\times P_5\times I^{n-4}$.
Since  $\beta^{-1, 2j}$ is equal to the number of degree $2j$ monomial
elements in a minimal basis of the Stanley-Reisner ideal of the polytope,
we can see that
$$\beta^{-1,2j}(P_5)=\begin{cases}5, & j=2\\ 0, & j\ge 3\end{cases}, \quad
\beta^{-1,2j}(P_6)=\begin{cases}9, & j=2\\ 0, & j\ge 3\end{cases},\quad$$
$$
\beta^{-1,2j}(Q)=\begin{cases}10, & j=2\\ 0, & j\ge
3\end{cases},\quad \beta^{-1,2j}(I^k)=\begin{cases}k, & j=2\\ 0, &
j\ge 3\end{cases}.$$ By \autoref{propositon: relations of Betti
numbers} (3), $\beta^{-1,2j}(P'\times I^k)=0$ for $j\ge 3$ where
$P'\approx I^2$, $P_5$, $P_6$, $Q$, or $P_5\times P_5$.
\end{proof}

We now prepare for the proof of \autoref{theorem: trianlgefree polytope}.
By \autoref{theorem: Hochster} we have
$\beta^{-2, 8}(P)=\sum_{|\sigma|=4}\dim\widetilde H^1(P_\sigma)$.
Therefore
\begin{eqnarray*}
\beta^{-2,8}(P_5) & = & \beta^{-2,8}(P_6) = 0,\\
\beta^{-2,8}(Q)   & =  &5,\\
\beta^{-2,8}(P_5\times P_5)&=&\beta^{-1,4}(P_5)\beta^{-1,4}(P_5)=25
\end{eqnarray*}
(note that since $Q$ does not have triangular faces,
$\beta^{-2,8}(Q)$ equals the number of 4-facet ``belts'' in~$Q$).
Hence we have
\begin{eqnarray*}
\beta^{-1,4}(P_6\times I^{n-2}) &= & \beta^{-1,4}(P_6)+\beta^{-1,4}(I^{n-2})=n+7,\\
\beta^{-2, 8}(P_6\times I^{n-2}) &= & \beta^{-1,4}(P_6)\cdot\beta^{-1,4}(I^{n-2})+
\beta^{0,0}(P_6)\cdot\beta^{-2,8}(I^{n-2})\\
& & +\,\beta^{-2,8}(P_6)\cdot\beta^{0,0}(I^{n-2}).
\end{eqnarray*}
On the other hand, by an inductive application of
\autoref{propositon: relations of Betti numbers} (3) we can see
easily that $\beta^{-2,8}(I^{n-2})=(n-2)(n-3)/2$. Therefore we
have
$$\beta^{-1,4}(P_6\times I^{n-2})=n+7, \quad \beta^{-2,8}(P_6\times I^{n-2})=\frac{n^2+13n-30}2.$$
By a similar computation we have
\begin{align*}
&\beta^{-1,4}(Q\times I^{n-3})=n+7, && \beta^{-2,8}(Q\times I^{n-3})=\frac{n^2+13n-38}2\\
&\beta^{-1,4}(P_5\times P_5\times I^{n-4})=n+6, &&
\beta^{-2,8}(P_5\times P_5\times I^{n-4})=\frac{n^2+11n-10}2.
\end{align*}

\begin{theorem}
Every triangle-free $n$-dimensional simple convex polytope with less than $2n+3$ facets
is rigid.
\end{theorem}
\begin{proof}
Let $P$ be  triangle-free with $f_0(P)\le 2n+2$, and let $M$ be a
quasitoric manifold over $P$. Let $P'$ be another simple convex
polytope and $M'$ a quasitoric manifold over $P'$. If
$H^\ast(M:\mathbb Q)\cong
H^\ast(M':\mathbb Q)$ as graded rings, then by
\autoref{proposition: bigraded Betti number} we have the equality
$\beta^{-i, 2j}(P)=\beta^{-i,2j}(P')$ for all $1\le i, j\le m$.
Since $P$ is triangle-free with $f_0(P)\le 2n+2$, $\beta^{-1,
2j}(P)=0$ for all $j\ge 3$ by \autoref{lemma: necessatity for
triangle-free}. Hence $\beta^{-1, 2j}(P')=0$ for all $j\ge 3$, and
\autoref{lemma:minimal basis and triangle-free} implies that $P'$
is triangle-free. Furthermore $H^\ast(M:\mathbb Q)\cong
H^\ast(M':\mathbb Q)$ implies in particular $f_0(P)=f_0(P')$. If
$f_0(P)=2n$ or $2n+1$, then there is only one simple
polytope with the given number of facets.
So $P\approx P'$. When
$f_0(P)=2n+2$ then there are three possible polytopes, but the
above computation shows that $\beta^{-i, 2j}$ are distinct for
these three polytopes. This shows that $P\approx P'$, which proves
the theorem.
The existence of quasitoric manifolds over $P$ is clear because we know
the existence of quasitoric manifolds over any two or three dimensional
simple convex polytopes and any $n$-simplex as well as any finite product of these
polytopes.
\end{proof}

\section{Rigidity of products of simplices}
We will make use of the following invariant in this  and the next
section.
\begin{definition}
The \emph{sigma invariant} of $P$ is $\sigma(P)=\sum_{j\ge
2}j\beta^{-1,2j}(P)$.
\end{definition}
\autoref{proposition: bigraded Betti number} implies that
$\sigma(P)$ is a cohomology invariant of quasitoric manifolds over
$P$. As we observed in \autoref{section: bigraded Betti number}
the Betti number $\beta^{-1,2j}(P)$ is equal to the number of
degree $2j$ elements in a minimal basis of the Stanley-Reisner
ideal $I_P$ of $P$. Therefore $2 \sigma(P)$ is nothing but the sum
of the degrees of all elements of a minimal basis of $I_P$.

\begin{lemma}\label{lemma:product of simplices}
Let $P$ be a simple polytope with $m$ facets. Then the following
conditions are equivalent:
\begin{itemize}
\item[(a)] $\sigma(P)=m$;
\item[(b)] the canonical minimal basis $\mathcal B$ of $I_P$ forms
a regular sequence;
\item[(c)] $P$ is combinatorially equivalent to a product of
simplices.
\end{itemize}
\end{lemma}

\begin{proof}
$(\mathrm c)\Rightarrow(\mathrm a)$ Clear.

$(\mathrm a)\Rightarrow(\mathrm b)$ Let $\mathbb Q(P)=\mathbb
Q[x_1,\ldots, x_m]/I_P$, where $x_i$ corresponds to a facet $F_i$
of $P$. Let $\mathcal B=\{g_1,\ldots, g_t\}$ be the canonical  minimal basis of
$I_P$. Since $\sigma(\mathbb Q(P))=m$, each $x_j$ must appear in
exactly one element of $\mathcal B$ with exponent $1$. It follows
easily that $g_1,\ldots,g_t$ is a regular sequence.

$(\mathrm b)\Rightarrow(\mathrm c)$. Let $\mathbb Q(P)=\mathbb
Q[x_1,\ldots, x_m]/(g_1,\ldots, g_t)$, where $g_1,\ldots,g_t$ is a
monomial regular sequence. It is well known~\cite[\S3.2]{BP} that
$g_1,\ldots,g_t$ is a regular sequence if and only if $g_i$ is not
a zero divisor in the quotient $\mathbb
Q[x_1,\ldots,x_m]/(g_1,\ldots,g_{i-1})$ for $1\le i\le t$ (this
property is often taken as the definition of a regular sequence).
Assume that some $x_j$ appears in more than one of
$g_1,\ldots,g_t$, say in $g_1$ and $g_2$. Then $g_2$ is a zero
divisor in $\mathbb Q[x_1,\ldots,x_m]/(g_1)$, which leads to a
contradiction. Therefore, each $x_j$ appears in at most one of the
monomials $g_1,\ldots,g_t$. Since every $x_j$ must appear in at
least one element in $I_P$, we obtain that every $x_j$ enters in
exactly one of $g_1,\ldots,g_t$. So we can rename $x_1,\ldots,x_m$
by
$y_{1\,0},\ldots,y_{1\,n_{\!1}},\ldots,y_{t\,0}\ldots,y_{t\,n_t}$
such that $g_j=\prod_{k=0}^{n_j}y_{j\,k}$  for $j=1,\ldots, t$.
Therefore we can see immediately that
\begin{align*}
\mathbb Q(P)&\cong \mathbb Q[x_1,\ldots, x_m]/I_P\\
&\cong \otimes_{i=1}^t\mathbb Q[y_{i\,0},\ldots,y_{i\,{n_i}}]/(g_i)\\
&\cong \otimes_{i=1}^t \mathbb Q(\Delta^{n_i})\\
&\cong \mathbb Q\Bigl(\prod_{i=1}^t\Delta^{n_i}\Bigr).
\end{align*}
Since the Stanley-Reisner ring with $\mathbb Q$-coefficients determines the combinatorial type of a
simple polytope \cite{BG}, we have $P\approx \prod_{i=1}^t\Delta^{n_i}$.
\end{proof}

Note that (b) in Lemma~\ref{lemma:product of simplices} is equivalent to saying
that $\mathbb Q(P)$ is a complete intersection ring.

\begin{theorem}
A finite product of simplices is rigid.
\end{theorem}
\begin{proof}
Let $M$ be a $2n$-dimensional quasitoric manifold over $P=\prod_{i=1}^t\Delta^{n_i}$.
Let $N$ be an another quasitoric manifold over a simple convex polytope $Q$, such that
$H^\ast(M:\mathbb Z)\cong H^\ast(N:\mathbb Z)$. Then
$H^\ast(M:\mathbb Q)\cong H^\ast(N:\mathbb Q)$ and
$f_i(P)=f_i(Q)$ for all $i$. In particular, $\sigma(\mathbb Q(P))=f_{0}(P)=f_{0}(Q)=n+t$.
Thus $Q$ is a simple convex polytope with $\sigma(\mathbb Q(Q))=f_{0}(Q)$.
Therefore $Q$ is also a product of simplices, i.e.,
$Q\approx \prod_{j=1}^s\Delta^{m_j}$.
But $H^\ast(M:\mathbb Q)\cong H^\ast(N:\mathbb Q)$ implies $\beta^{-1,2j}(P)=\beta^{-1,2j}(Q)$
for all $J$. This implies that $\{n_i\}=\{m_j\}$ and $t=s$. Thus $P\cong Q$.
\end{proof}

\section{Rigidity of vertex cuts}

The following proposition shows that certain Betti numbers and the
sigma invariant of a vertex cut of $P$ are independent of the
choice of the cut vertex; whereas the combinatorial type of $P$
may depend of this choice, see Example~\ref{vcex}.

\begin{proposition}\label{proposition: sigma and Betti number of vertex cut}
Let $P$ be an $n$-dimensional simple convex polytope with $m$
facets, which is different from the $n$-simplex $\Delta^n$.
Then we have\\[-6pt]
\begin{enumerate}
\item $\beta^{-1, 2j}(\vc(P))=
    \begin{cases} \beta^{-1, 2j}(P) + m-n, & j=2\\
                  \beta^{-1,2j}(P),        & 3\le j\le n-1\\
                  \beta^{-1,2j}(P)+1,        & j=n
    \end{cases}$\\
\item $\sigma(\vc(P))=\sigma(P)+2m-n$.
\end{enumerate}
\end{proposition}
\begin{proof}
Both statements follow easily from the interpretation
of~$\beta^{-1, 2j}(P)$ as the number of degree $2j$ elements in
the minimal basis of the ideal~$I_P$.
\end{proof}

When $P=\prod_{i=1}^t\Delta^{n_i}$ with $t\ne 1$, we have $n=\sum_{i=1}^t n_i$, $m=n+t$ and
$\sigma(P)=m$. Hence we have $\sigma(\vc(P))=3m-n$.

Let $\mathcal F=\{F_1,\ldots, F_m\}$ be the set of facets in $P$.
Let $x_i$ be the corresponding generator to $F_i$ in $\mathbb
Q(P)=\mathbb Q[x_1,\ldots, x_m]/I_P$. Let ${\mathcal
B}=\{h_1,\ldots, h_\ell\}$ be the canonical minimal basis for
$I_P$. For each $x_i$ the {\em frequency} $\mathfrak{f}(x_i)$ is
the number of $h_k$ in ${\mathcal B}$ divisible by $x_i$.

\begin{lemma}\label{lemma: decomposable polytopes}
Let $P$ be an $n$-dimensional simple convex polytope. Let
${\mathcal B}$ be the canonical minimal basis for $I_P$. If
$\mathfrak{f}(x_i)=1$ for some $i$, then $P \approx\Delta^k\times
P' $ for some polytope $P'$ of dimension $n-k$ and $k=\deg h/2 -1$
where $h$ is the unique element in ${\mathcal B}$ such that
$x_i|h$.
\end{lemma}
\begin{proof}
Let $\mathcal{B}=\{h_1, \ldots, h_s \}$. Assume
$\mathfrak{f}(x_1)=1$ and  $h_1 = x_1 \cdots x_t$ for simplicity.
Hence $h_1$ is the unique element of ${\mathcal B}$ that is
divisible by $x_1$. We claim that $\mathfrak{f}(x_2) = \cdots =
\mathfrak{f}(x_t)=1$. Assume otherwise, say $\mathfrak{f}(x_2)
\geq 2$. Without loss of generality, we may assume $h_2 = x_2
x_{i_1} \cdots x_{i_k}$. Then $x_{i_j}\ne x_1$ and $x_2$ for all
$j=1,\ldots, k$ because $h_1$ is the only element of ${\mathcal
B}$ divisible by $x_1$. Since $h_2\in {\mathcal B}$, if we let
$T:=F_{i_1} \cap \cdots \cap F_{i_k}$ then $T\ne \emptyset$ but
$F_2\cap T=\emptyset$. On the other hand, since $x_1\nmid h_2$ we
have $F_1\cap T\ne \emptyset$
(otherwise $x_1x_{i_1}\cdots x_{i_k}\in I_P$, so there would be
another element in $\mathcal B$ divisible by~$x_1$).

If $k\ge n$, then $F_1\cap T\ne \emptyset$ implies that more than
$n$ facets of $P$ are intersecting, which is impossible because
$P$ is simple. Therefore $\dim T=n-k\ge 1$. Since $\dim(F_1\cap
T)=\dim T-1$, there exists a
vertex $v$ of $T$ which does
not belong to $F_1$. Let $v$ be the intersection of $n$ facets
$F_{\ell_1},\ldots, F_{\ell_n}$. Since $F_2\cap T=\emptyset$, the
vertex $v$ does not belong to $F_2$, hence $F_{\ell_j}\ne F_2$ for
all $j=1,\ldots, n$. Since $v$ does not belong to $F_1$, we have
$F_1\cap F_{\ell_1}\cap \cdots\cap F_{\ell_n}=\emptyset$.
Therefore there must exist an element $h\in{\mathcal B}$, which
divides the monomial $x_1x_{\ell_1}\cdots x_{\ell_n}$. But since
$F_{\ell_1}\cap\cdots\cap F_{\ell_n}=v\ne\emptyset$, the element
$h$ must be divisible by $x_1$. Since $F_{\ell_j}\ne F_2$ for all
$j=1,\ldots, n$, it follows that $x_2\nmid h$. Thus $h_1\nmid h$,
which contradicts to the condition that $\mathfrak f(x_1)=1$. This
shows that $\mathfrak f(x_2)=1$, and by a similar argument we can
see that $\mathfrak f(x_i)=1$ for all $i=1,\ldots, t$. Hence,
$$
    \Q(P) = \Q[x_1, \cdots, x_t]/h_1 \otimes \Q[x_{t+1}, \cdots, x_m]/I',
$$ where $I'$ is the ideal generated by $\mathcal{B} \setminus \{ h_1 \}$.

Since $\mathbb Q[x_1,\ldots x_t]/(h_1)\cong\mathbb Q(\Delta^t)$,
it is enough to prove that there is an isomorphism $\mathbb
Q[x_{t+1},\ldots, t_m]/I'\cong \mathbb Q(P')$ for some polytope
$P'$ of dimension $n-k$. (Indeed, then we instantly get
$P\approx\Delta^k\times P'$ because
the
rational Stanley-Reisner ring determines the combinatorial
type of a simple polytope. \cite{BG}) Let $P' := F_2 \cap \cdots \cap
F_{t}$. Then every facet except $F_1$ intersects with $P'$. Let
$G_j=F_j\cap P'$ for $j=t+1, \ldots, m$. Then $G_j$'s are facets
of $P'$. This implies that the face poset structure of $P'$ agrees
the face poset structure of $\{F_{t+1} , \ldots , F_{m} \}$.  Thus
$\mathcal{B} \setminus \{ h_1 \} = \{ h_2 , \ldots, h_s \}$ is the
canonical minimal basis for $I_{P'}$. Hence $\mathbb
Q[x_{t-1},\ldots, x_n]/I'\cong \mathbb Q(P')$.
\end{proof}

\begin{theorem}\label{theorem: characterization of cv(prod of simp)}
Let $Q$ be an $n$-dimensional simple convex polytope with $m+1$ facets. If $\sigma(Q)=3m-n$ and $\beta^{-1, 2n}(Q)\neq 0$, then $Q$ is a vertex cut of a product of simplices.
\end{theorem}
\begin{proof}
We claim that one of the facets of $Q$ is an $(n-1)$-simplex. Then
$Q$ is a vertex cut of some simple convex polytope $P$. By
\autoref{proposition: sigma and Betti number of vertex cut} we
have
$$
    \sigma(P) = \sigma(Q) - 2m+n = (3m-n)-(2m-n) = m.
$$
Thus by \autoref{lemma:product of simplices} $P$ is a product of
simplices, and we are done. We now prove the claim. Let
$F_1$,\ldots, $F_{m+1}$ be the facets of $Q$ and let $x_1$,
\ldots, $x_{m+1}$ be the associated generators of $\mathbb Q(Q)$.
Let ${\mathcal B}$ be the canonical minimal basis for the ideal
$I_Q$. Since $\beta^{-1, 2n}(Q) \geq 1$, there exists $\widetilde
h\in{\mathcal B}$ with $\deg \widetilde h=2n$. Without loss of
generality we may assume $\widetilde h=x_1\ldots x_n$. Then we can
see easily that $F_1\cup\cdots\cup F_n$ is homeomorphic to
$S^{n-2}\times I$, while $F_1\cup\cdots\cup F_{m+1}\cong S^{n-1}$.
Thus $F_1\cup\cdots\cup F_{m+1}\backslash F_1\cup\cdots\cup
F_n=F_{n+1}\cup\cdots \cup F_{m+1}$ has two connected components.
For simplicity let $F_{n+1}\cup\cdots \cup F_{n+k}$ and
$F_{n+k+1}\cup\cdots \cup F_{m+1}$ be the two components. Then
$F_{n+i}\cap F_{n+j}=\emptyset$ for $i=1,\ldots, k$ and
$j=k+1,\ldots, m+1-n$.

If $k=1$ or $m-n$, then one of the components of
$F_{n+1}\cup\cdots\cup F_{m+1}$ is a
single
facet of $Q$,
and this facet is an $(n-1)$-simplex. This proves the claim.
Assume otherwise, i.e., suppose $2\le k\le [(m+1-n)/2]$. Let
$\mathcal B_1=\{x_{n+i}x_{n+j} | i=1,\ldots, k \text{ and
}j=k+1,\ldots, m+1-n \}$. Then we have
\begin{equation}\label{equation: inequality for B1}
\sum_{h\in \mathcal B_1}\deg(h) = 4k(m+1-n-k) \geq 8(m-n-1)
\end{equation}
because $k(m+1-n-k)$ is increasing for $2\le k\le [(m+1-n)/2]$. Note
that
the frequencies satisfy $\mathfrak{f}(x_i) \geq 2$ for all $i=1,
\ldots, n$ because otherwise \autoref{lemma: decomposable
polytopes}
would imply that $Q \approx \Delta^{n-1} \times \Delta^{1}$, but
in this case $\sigma(Q) =n+2 \neq 3m-n = 2n+3$. Therefore for each
$x_i$, there
exists $h_i \in \mathcal{B}$ such that $x_i |h_i$ and $h_i \neq
\widetilde{h}$ for $i=1,\ldots,n$. Note that some of
$h_1,\ldots,h_n$ may coincide. So we let $\widetilde{h}_1, \ldots,
\widetilde h_s$ denote  all distinct elements among $h_i$'s. If
$s=1$, then $\widetilde h_1$ is divisible by all $x_i$ for $i=1,
\ldots, n$. Hence $\widetilde{h} | \widetilde h_1$ and therefore
$\widetilde{h}=\widetilde h_1$, which is a contradiction.
Therefore $s \geq 2$.

If $s \geq 3$, then
\begin{eqnarray} \label{equation: inequality for B2_1}
\sum_{h\in \mathcal B \setminus \mathcal B_1} \deg h & \ge & \deg
\widetilde h + \sum_{i=1}^s\deg \widetilde h_i\\\nonumber &\ge &
\deg\widetilde h + 2n + 6= 4n + 6,
\end{eqnarray}
where the last inequality follows from the conditions $s\ge 3$,
$\deg \widetilde h_i\ge 4$ and $x_1\cdots x_n\mid \widetilde
h_1\cdots\widetilde h_s$.

Suppose $s = 2$. Then without loss of generality we may assume
that $\widetilde h_1 = g_1 x_1 \cdots x_\ell$ and $\widetilde
h_2=g_2 x_{\ell+1} \cdots x_n$  with
$1\le \ell \le n-1$ for some monomials $g_1, g_2$ in $x_{n+1},
\ldots, x_{m+1}$ of degree $\geq 2$. If degree $g_2 \geq 4$, then
$$\deg \widetilde h_1 + \deg
\widetilde h_2 = 2n + \deg g_1 + \deg g_2 \geq 2n+6.$$ Therefore
the inequality (\ref{equation: inequality for B2_1}) holds in this
case. Now suppose degree $g_2 = 2$. Then $g_2=x_i$ for $n+1\le
i\le n+k$ or $g_2=x_{n+j}$ for $k+1\le j\le m+1-n$. We only prove
the case when $g_2 = x_{n+k+1}$. The other cases are similar. In
this case consider the monomial $q = \prod_{j=2}^{n+2}x_j$. By the
assumption, $\widetilde h_i\nmid q$ for $i=1,2$.
But $q$ must vanish in $\Q(Q)$ because any set of $n+1$ facets has
empty intersection in a simple polytope.
Therefore there exists a monomial $q'$ of degree $\geq 4$ in
${\mathcal B} \setminus \mathcal B_1$, which divides $q$. Thus
$$\sum_{h \in
{\mathcal B} \setminus \mathcal B_1}\deg h\ge \deg \widetilde{h} +
\deg q' + \deg \widetilde h_1 + \deg\widetilde h_2 \geq 4n+8
>4n+6.$$

We thus have proved that in all cases
\begin{equation}\label{equation: inequality for both}
 2\sigma(Q) = \sum_{h\in \mathcal B} \deg h \ge 8(m-n-1) + 4n+6 = 8m - 4n - 2 .
\end{equation}

On the other hand,
by the assumption of the theorem, $\sigma(Q) = 3m-n$. Thus $3m-n
\geq 4m-2n-1$, hence $n+2 \ge m+1$. Therefore $Q$ is
combinatorially equivalent to either $\Delta^{n_1} \times
\Delta^{n_2}$ or
$\Delta^n$. But $\beta^{-1,2n}(Q)
\neq 0$
gives that $Q \approx \Delta^{n-1} \times \Delta^1$, which implies
 $\sigma(Q) = m+1 \neq 3m-n$. This is a
contradiction. Thus we have $k=1$ or $m-n$, which proves the
theorem.
\end{proof}

\begin{theorem}
If
$P$ is a finite product of simplices, then $\vc(P)$ is  rigid.
\end{theorem}

\begin{proof}
If $P$ is an $n$-simplex, then $\vc(P) = \Delta^{n-1}\times
\Delta^1$, which is rigid by \autoref{theorem: product of
simplices is rigid}. Assume otherwise. Let $Q=\vc(P)$, and let $M$
be a quasitoric manifold over $Q$. Suppose $N$ is quasitoric
manifold over another simple convex polytope $Q'$ such that
$H^\ast(M:\mathbb Q)\cong H^\ast(N:\mathbb Q)$ as graded rings.
Then $\beta^{-1, 2j}(Q)=\beta^{-1,2j}(Q')$, and hence
$\sigma(Q')=\sigma(Q)=3m-n$ and
$\beta^{-1,2n}(Q')=\beta^{-1,2n}(Q) \neq 0$. By \autoref{theorem:
characterization of cv(prod of simp)} $Q'=\vc(P')$ for $P' = \prod
\Delta^{n_i}$. By \autoref{proposition: sigma and Betti number of
vertex cut}(1), $\beta^{-1, 2j}(Q)=\beta^{-1,2j}(Q')$ implies
$\beta^{-1,2j}(P) = \beta^{-1,2j}(P')$ for all $j$. Both $P$ and
$P'$ are products of simplices, thus $P \approx P'$. So $Q \approx
Q'$, which proves the theorem.
\end{proof}

\section{Rigidity of $3$-dimensional simple convex polytopes}\label{section:3-polytope}

Since rigidity of $2$-dimensional simple convex polytope is
settled by Corollary~\ref{corollary: polygon is rigid}, rigidity
of $3$-dimensional simple convex polytope is naturally the next
target. Note that any $3$-dimensional simple convex polytope
supports a quasitoric manifolds. The four color problem gives an
easy proof of this.

\psset{unit=4pt}
\begin{table}
  \centering
\begin{tabular}{|c|c|c|}
  \hline
\begin{tabular}{c} Type \\ \tiny{(%
Betti numbers)} \\ \end{tabular} & Simple polytopes &  \\ \hline
\hline
\begin{tabular}{c} $3^4$ \\ \tiny{()} \\ \end{tabular} & \begin{pspicture}(0,-1)(6,6)
\pspolygon(3,0)(0.4,4.5)(5.6,4.5) \psline(0.4,4.5)(3,3)
\psline(3,0)(3,3) \psline(5.6,4.5)(3,3) \rput(3,-1){\tiny{$3^4$}}
\end{pspicture} & rigid \\ \hline
\begin{tabular}{c} $\vc(3^4)$ \\ \tiny{(1)} \\ \end{tabular} & \begin{pspicture}(0,-1)(6,6)
\pspolygon(3,0)(0.4,4.5)(5.6,4.5) \pspolygon(3,1.5)(4.3,3.75)(1.7,
3.75) \psline(0.4,4.5)(1.7,3.75) \psline(3,0)(3,1.5)
\psline(5.6,4.5)(4.3,3.75) \rput(3,-1){\tiny{$3^2 4^3$}}
\end{pspicture} & rigid \\ \hline
\begin{tabular}{c} $\vc^2(3^4)$ \\ \tiny{(3,2)} \\ \end{tabular} & \begin{pspicture}(0,-1)(6,6)
\pspolygon(3,6)(0.15,3.93)(1.24,0.57)(4.76,0.57)(5.85,3.93)
\pspolygon(3,4.5)(1.57,3.46)(1.24,0.57)(4.76,0.57)(4.43,3.46)
\psline(3,6)(3,4.5) \psline(0.15,3.93)(1.57,3.46)
\psline(5.85,3.93)(4.43,3.46) \rput(3,-1){\tiny{$3^24^25^2$}}
\end{pspicture} & rigid \\ \hline
\begin{tabular}{c} $4^6$ \\ \tiny{(3,0)} \\ \end{tabular} & \begin{pspicture}(0,-1)(6,6)
\pspolygon(5.12,5.12)(0.88,5.12)(0.88,0.88)(5.12,0.88)
\pspolygon(4.06,4.06)(1.94,4.06)(1.94,1.94)(4.06,1.94)
\psline(5.12,5.12)(4.06,4.06) \psline(0.88,5.12)(1.94,4.06)
\psline(0.88,0.88)(1.94,1.94) \psline(5.12,0.88)(4.06,1.94)
\rput(3,-1){\tiny{$4^6$}}
\end{pspicture} & rigid \\ \hline
\begin{tabular}{c} $\vc^3(3^4) $ \\ \tiny{(6,8,3)} \\ \end{tabular} & $3^2 4^3 6^2$, $3^35^36^1$,$3^24^25^26^1$ (3 polytopes) & nonrigid \\ \hline
\begin{tabular}{c} $\vc(4^6) $ \\ \tiny{(6,6,1)} \\ \end{tabular} & \begin{pspicture}(0,-1)(6,6)
\pspolygon(3,6)(0.15,3.93)(1.24,0.57)(4.76,0.57)(5.85,3.93)
\pspolygon(3,4.5)(1.57,3.46)(3,2.09)(4.43,3.46)
\psline(1.24,0.57)(3,1.49) \psline(4.76,0.57)(3,1.49)
\psline(3,2.09)(3,1.49) \psline(3,4.5)(3,6)
\psline(0.15,3.93)(1.57,3.46)\psline(5.85,3.93)(4.43,3.46)
\rput(3,-1){\tiny{$3^1 4^3 5^3$}}
\end{pspicture} & rigid \\ \hline

\begin{tabular}{c} $4^5 5^2$ \\ \tiny{(6,5,0)} \\ \end{tabular} & \begin{pspicture}(0,-1)(6,6)
\pspolygon(3,6)(0.15,3.93)(1.24,0.57)(4.76,0.57)(5.85,3.93)
\pspolygon(3,4.5)(1.57,3.46)(2.12,1.79)(3.88,1.79)(4.43,3.46)
\psline(3,6)(3,4.5) \psline(0.15,3.93)(1.57,3.46)
\psline(1.24,0.57)(2.12,1.79) \psline(4.76,0.57)(3.88,1.79)
\psline(5.85,3.93)(4.43,3.46) \rput(3,-1){\tiny{$4^55^2$}}
\end{pspicture} & rigid \\ \hline

\begin{tabular}{c} $\vc^4(3^4)$ \\ \tiny{(10,20,15,4)} \\ \end{tabular} & \small{\begin{tabular}{l} $3^2 4^4 7^2$, $3^34^15^26^17^1$, $3^24^35^16^17^1$, $3^24^25^37^1$, $3^46^4$ \\ $3^34^15^16^3$, $3^24^25^26^2(i)$ (7 polytopes) \end{tabular}}  & nonrigid \\ \hline

\begin{tabular}{c} $\vc^2(4^6) $ \\ \tiny{(10,18,11,2)} \\ \end{tabular} & \small{$3^2 4^2 5^2 6^2(ii)$, $3^1 4^45^16^2$, $3^24^15^46^1$, $3^2 5^6$ (4 polytopes)} & nonrigid \\ \hline

\begin{tabular}{c} $\vc(4^5 5^2)$ \\ \tiny{(10,17,9,1)} \\ \end{tabular} & \begin{pspicture}(0,-1)(6,6)
\pspolygon(0,3)(1.5,0.4)(4.5,0.4)(6,3)(4.5,5.6)(1.5,5.6)
\pspolygon(1.5,3)(3,2)(4.5,3)(3.75,4.3)(2.25,4.3)
\psline(0,3)(1.5,3) \psline(1.5,0.4)(3,1.39)
\psline(4.5,0.4)(3,1.39) \psline(6,3)(4.5,3)
\psline(4.5,5.6)(3.75,4.3) \psline(1.5,5.6)(2.25,4.3)
\psline(3,1.39)(3,2) \rput(3,-1){\tiny{$3^1 4^3 5^3 6^1$}}
\end{pspicture} & rigid \\ \hline

\begin{tabular}{c} $4^6 6^2$  \\ \tiny{(10,16,9,0)} \\ \end{tabular} & \begin{pspicture}(0,-1)(6,6)
\pspolygon(0,3)(1.5,0.4)(4.5,0.4)(6,3)(4.5,5.6)(1.5,5.6)
\pspolygon(1.5,3)(2.25,1.7)(3.75,1.7)(4.5,3)(3.75,4.3)(2.25,4.3)
\psline(0,3)(1.5,3) \psline(1.5,0.4)(2.25,1.7)
\psline(4.5,0.4)(3.75,1.7) \psline(6,3)(4.5,3)
\psline(4.5,5.6)(3.75,4.3) \psline(1.5,5.6)(2.25,4.3)
\rput(3,-1){\tiny{$4^66^2$}}
\end{pspicture} & rigid \\ \hline

\begin{tabular}{c} $4^4 5^4$  \\ \tiny{(10,16,5,0)} \\ \end{tabular} & \begin{pspicture}(0,-1)(6,6)
\pspolygon(3,6)(0.15,3.93)(1.24,0.57)(4.76,0.57)(5.85,3.93)
\pspolygon(3,4.5)(1.57,3.46)(2.12,2.09)(3.88,2.09)(4.43,3.46)
\pspolygon(2.12,2.09)(3.88,2.09)(3.88,1.49)(2.12,1.49)
\psline(3,6)(3,4.5) \psline(0.15,3.93)(1.57,3.46)
\psline(1.24,0.57)(2.12,1.49) \psline(4.76,0.57)(3.88,1.49)
\psline(5.85,3.93)(4.43,3.46) \rput(3,-1){\tiny{$4^45^4$}}
\end{pspicture} & rigid \\ \hline
\end{tabular}
\
  \caption{Rigidity of simple 3-polytopes with $f_0 \leq 8$} \label{table:f_0 <= 8}
\end{table}

\begin{table}
  \centering

\begin{tabular}{|c|c|c|}
  \hline
\begin{tabular}{c} Type \\ \tiny{(%
Betti numbers)} \\ \end{tabular} & Simple polytopes &  \\ \hline
\hline
\begin{tabular}{c} $\vc^5(3^4)$ \\ \tiny{(15,40,45,24,5)} \\ \end{tabular}  & \tiny{\begin{tabular}{c} $3^24^58^2$, $3^34^25^27^18^1(i)$, $3^34^25^27^18^1(ii)$, $3^24^25^17^18^1$, $3^45^26^28^1$\\
$3^34^25^16^28^1(i)$, $3^34^25^16^28^1(ii)$, $3^24^46^28^1$, $3^34^15^36^18^1$, $3^24^35^26^18^1(i)$\\
$3^24^35^26^18^1(ii)$, $3^24^25^48^1$, $3^44^16^27^2$, $3^45^26^17^2$, $3^34^25^16^17^2(i)$\\
$3^34^25^16^17^2(ii)$, $3^34^15^37^2$, $3^24^35^27^2(ii)$, $3^34^26^37^1(i)$, $3^34^26^37^1(ii)$\\
$3^34^15^26^27^1(ii)$, $ 3^24^35^16^27^1(iii)$,
$3^24^25^36^17^1(iii)$, $3^24^25^26^3(iii)$ \\(24 polytopes)
\end{tabular}
 }    & nonrigid \\ \hline
\begin{tabular}{c} $\vc^3(4^6)$ \\ \tiny{(15,38,39,18,3)} \\ \end{tabular}  &\tiny{\begin{tabular}{c} $3^24^35^27^2(i)$, $3^14^55^17^2$, $3^34^15^26^27^1(i)$, $3^24^35^16^27^1(i)$, $3^24^35^16^27^1(ii)$ \\ $3^24^25^36^17^1(ii)$, $3^14^45^26^17^1(iii)$, $3^35^36^3(i)$, $3^35^36^3(ii)$, $3^24^25^26^3(ii)$ \\ $3^24^15^46^2(ii)$ (11 polytopes)\end{tabular}}& nonrigid \\ \hline
\begin{tabular}{c} $\vc^2(4^55^2)$ \\ \tiny{(15,37,36,15,2)} \\ \end{tabular}  & \tiny{\begin{tabular}{c} $3^24^25^36^17^1(i)$, $3^14^45^26^17^1(i)$, $3^24^15^57^1$, $3^24^36^4$, $3^24^25^26^3(i)$ \\ $3^14^45^16^3$, $3^24^15^46^2(i)$ (7 polytopes)\end{tabular}} & nonrigid \\ \hline
\begin{tabular}{c} $\vc(4^66^2) $\\ \tiny{(15,36,35,14,1)} \\ \end{tabular} & \begin{pspicture}(0,-1)(6,6)
\pspolygon(3,6)(0.65,4.87)(0.08,2.33)(1.69,0.30)(4.31,0.30)(5.92,2.33)(5.35,4.87)
\pspolygon(3,4.5)(1.83,3.93)(1.54,2.67)(3,1.95)(4.46,2.67)(4.17,3.93)
\psline(3,6)(3,4.5) \psline(0.65,4.87)(1.83,3.93)
\psline(0.08,2.33)(1.54,2.67) \psline(1.69,0.30)(3,1.35)
\psline(4.31,0.30)(3,1.35) \psline(5.92,2.33)(4.46,2.67)
\psline(5.35,4.87)(4.17,3.93) \psline(3,1.35)(3,1.95)
\rput(3,-1){\tiny{$3^1 4^4 5^2 6^1 7^1 (ii)$}}
\end{pspicture} & rigid \\ \hline
\begin{tabular}{c} $\vc(4^4 5^4)$ \\ \tiny{(15,36,31,10,1)} \\ \end{tabular}   & $3^14^35^36^2$, $3^14^25^56^1$ (2 polytopes) & nonrigid \\ \hline
\begin{tabular}{c} $4^77^2$ \\ \tiny{(15,35,35,14,0)} \\ \end{tabular}  & \begin{pspicture}(0,-1)(6,6)
\pspolygon(3,6)(0.65,4.87)(0.08,2.33)(1.69,0.30)(4.31,0.30)(5.92,2.33)(5.35,4.87)
\pspolygon(3,4.5)(1.83,3.93)(1.54,2.67)(2.35,1.65)(3.65,1.65)(4.46,2.67)(4.17,3.93)
\psline(3,6)(3,4.5) \psline(0.65,4.87)(1.83,3.93)
\psline(0.08,2.33)(1.54,2.67) \psline(1.69,0.30)(2.35,1.65)
\psline(4.31,0.30)(3.65,1.65) \psline(5.92,2.33)(4.46,2.67)
\psline(5.35,4.87)(4.17,3.93) \rput(3,-1){\tiny{$4^7 7^2$}}
\end{pspicture} & rigid \\ \hline
\begin{tabular}{c} $4^6 \sharp 4^6$ \\ \tiny{(15,36,33,12,1)} \\ \end{tabular} & \begin{pspicture}(0,-1)(6,6)
\pspolygon(0,3)(1.5,0.4)(4.5,0.4)(6,3)(4.5,5.6)(1.5,5.6)
\pspolygon(1,3)(1.5,3.87)(2,3)(1.5,2.13)
\pspolygon(4,3)(4.5,3.87)(5,3)(4.5,2.13) \psline(2,3)(4,3)
\psline(0,3)(1,3) \psline(5,3)(6,3) \psline(1.5,0.4)(1.5,2.13)
\psline(1.5,5.6)(1.5,3.87) \psline(4.5,0.4)(4.5,2.13)
\psline(4.5,5.6)(4.5,3.87) \rput(3,-1){\tiny{$4^6 6^3$}}
\end{pspicture} & rigid \\ \hline
\begin{tabular}{c} $4^5 5^2 6^2$ \\ \tiny{(15,35,29,8,0)} \\ \end{tabular}   & \begin{pspicture}(0,-1)(6,6)
\pspolygon(0,3)(1.5,0.4)(4.5,0.4)(6,3)(4.5,5.6)(1.5,5.6)
\pspolygon(1.5,3)(2.25,2)(3.75,2)(4.5,3)(3.75,4.3)(2.25,4.3)
\pspolygon(2.25,1.4)(3.75,1.4)(3.75,2)(2.25,2) \psline(0,3)(1.5,3)
\psline(1.5,0.4)(2.25,1.4) \psline(4.5,0.4)(3.75,1.4)
\psline(6,3)(4.5,3) \psline(4.5,5.6)(3.75,4.3)
\psline(1.5,5.6)(2.25,4.3) \rput(3,-1){\tiny{$4^5 5^2 6^2$}}
\end{pspicture} & rigid \\ \hline
\begin{tabular}{c} $4^4 5^4 6^1$ \\ \tiny{(15,35,27,6,0)} \\ \end{tabular}   & \begin{pspicture}(0,-1)(6,6)
\pspolygon(0,3)(1.5,0.4)(4.5,0.4)(6,3)(4.5,5.6)(1.5,5.6)
\pspolygon(1.5,3)(2.25,1.7)(3,2)(3.75,1.7)(4.5,3)(3.75,4.3)(3,4)(2.25,4.3)
\psline(3,2)(3,4) \psline(0,3)(1.5,3) \psline(1.5,0.4)(2.25,1.7)
\psline(4.5,0.4)(3.75,1.7) \psline(6,3)(4.5,3)
\psline(4.5,5.6)(3.75,4.3) \psline(1.5,5.6)(2.25,4.3)
\rput(3,-1){\tiny{$4^4 5^4 6^1$}}
\end{pspicture} & rigid\\ \hline
\begin{tabular}{c} $4^35^6$ \\ \tiny{(15,35,24,3,0)} \\ \end{tabular}   & \begin{pspicture}(0,-1)(6,6)
\pspolygon(3,6)(0.15,3.93)(1.24,0.57)(4.76,0.57)(5.85,3.93)
\pspolygon(3,4.5)(2.49,3.70)(1.57,3.46)(2.12,1.79)(3,2.24)(3.88,1.79)(4.43,3.46)(3.52,3.66)
\psline(3,6)(3,4.5) \psline(0.15,3.93)(1.57,3.46)
\psline(1.24,0.57)(2.12,1.79) \psline(4.76,0.57)(3.88,1.79)
\psline(5.85,3.93)(4.43,3.46) \psline(3,3)(3,2.24)
\psline(3,3)(2.49,3.70) \psline(3,3)(3.52,3.66)
\rput(3,-1){\tiny{$4^3 5^6$}}
\end{pspicture} & rigid \\ \hline
\end{tabular}
\
  \caption{Rigidity of simple 3-polytopes with $f_0 =9$}\label{table:f_0 = 9}
\end{table}


On pages 192 and 193 of Appendix A.5 in \cite{Od} there is a list
of 3-dimensional simple convex polytopes with $\le9$ facets. In
the list the polytopes are labeled in the form $\alpha^x \beta^y
\gamma^z$ which means the polytope has $x$ many $\alpha$-gon
facets, $y$ many $\beta$-gon facets, and $z$ many $\gamma$-gon
facets.
For example, the polytope $3^4$ is the tetrahedron, and $3^24^3$
is the triangular prism.

\autoref{table:f_0 <= 8} lists simple $3$-polytopes with $\le 8$
facets, their bigraded Betti numbers and rigidity.
\autoref{table:f_0 = 9} contains the same information about simple
3-polytopes with 9 facets.
In the tables
$\vc^k(P)$ denotes a $k$-fold vertex cut of~$P$. The Betti numbers
are listed in the form
\[
  (\beta^{-1,4},\ldots,\beta^{-(j-1),2j},\ldots,\beta^{-(m-4),2(m-3)}).
\]
Note that the numbers
above completely
determine all bigraded Betti numbers of a $3$-dimensional
polytope. Indeed, unless $(i,j)=(0,0)$ or $(m-3_, m)$ the number
$\beta^{-i, 2j}$ is zero for $j-i\ne 1,2$ by \autoref{theorem:
Hochster}. By \autoref{propositon: relations of Betti numbers}
$\beta^{0,0}=\beta^{-(m-3),2m}=1$ and
$\beta^{-(j-1),2j}=\beta^{-i', 2j'}$ where $i'=(m-3)-(j-1)$ and
$j'=m-j$. Note that $j' - i' = 2$. This implies
that the whole set of Betti numbers is determined
by $\beta^{-(j-1),2j}$'s for $j=2, \ldots, m$. Moreover for a
subset $\sigma \subset \{1, \ldots, m\}$, if $|\sigma| >m-3$, then
$P_\sigma$ is always connected.
We therefore consider only $\beta^{-(j-1),2j}$
for $j=2, \ldots,m-3$.

Both $3^4$ and $3^24^3$ are rigid by \autoref{theorem: product of
simplices is rigid}.

There are exactly two polytopes $3^24^25^2$ and $4^6$ with
$f_0=6$. Polytope $4^6$ is rigid because it is the cube $I^3$, and
polytope $3^24^25^2$ is the vertex cut of the triangular prism, so
it is also rigid by \autoref{theorem: cv of product of simplices}.

There are five different polytopes with $f_0=7$, which are
$3^24^36^2$, $3^35^36^1$, $3^24^25^26^1$, $3^14^35^3$, and $4^55^2$.
The first three are the polytopes obtained from  triangular prism
$\Delta^2\times I$ by taking vertex cuts twice. So by the  argument
of
Example~\ref{vcex} they are all nonrigid. The polytope $3^14^35^3$
is the vertex cut of the cube $I^3$ and hence rigid by
\autoref{theorem: cv of product of simplices}. Polytope $4^55^2$
is the pentagonal prism, which is rigid by \autoref{theorem:
trianlgefree polytope}.

There are $14$ different polytopes with $f_0=8$.  Seven of them
are obtained from the triangular prism by taking vertex cuts three
times, and so they are all nonrigid.
These
are $3^24^47^2$, $3^34^15^26^17^1$, $3^24^35^16^17^1$,
$3^24^25^37^1$, $3^46^4$, $3^34^15^16^3$, and $3^24^25^26^2(i)$.
There are four polytopes obtained from the cube by taking vertex
cuts twice. They are $3^24^25^26^2(ii)$, $3^14^55^16^2$,
$3^24^15^46^1$ and $3^25^6$, and all of them are nonrigid. The
remaining polytopes are $3^14^35^36^1$ which is the
vertex cut of the pentagonal prism, $4^66^2$ which is the
hexagonal prism $P_6\times I$, and $4^45^4$ which is obtained from
the pentagonal prism by cutting out a triangular prism shaped
neighborhood of an edge. Since the Betti numbers $\beta^{-1,2j}$
of $3^14^35^36^1$ are different from those of the other two and
also different from the previous groups, it is rigid. The
remaining polytopes are $4^66^2$ and $4^45^4$. These polytopes
have $2n+2=8$ facets. So by \autoref{theorem: trianlgefree
polytope} they are rigid.

There are $50$ different polytopes with $f_0 = 9$,
and only six of them are rigid.
Among them five are triangle-free polytopes, namely $4^66^3$,
$4^55^26^2$, $4^45^46^1$, $4^35^6$ and $4^77^2$, and the sixth is
the
polytope $3^14^45^26^17^1(ii)$ which is the vertex cut of $P_6
\times I$.
In each case the rigidity is established by comparing the Betti
numbers, and observing that these numbers for each of the six
polytopes are different from the other's.

Finally we give a proof of rigidity of a dodecahedron.
\begin{theorem}\label{theorem:dodecahedron}
    A dodecahedron is rigid.
\end{theorem}
\begin{proof}
A computation using Theorem~\ref{theorem: Hochster} shows that the
$(-2,8)$-th Betti number of a dodecahedron is 0. Let $P$ be a
simple 3-polytope with 12 facets whose Betti numbers are equal to
those of a dodecahedron. Let $x_k$ be the number of $k$-gonal
facet of $P$. By Euler equation $\sum_{k \geq 3} x_k (6-k) = 12$.
Since the number of facets $\sum_{k \geq 3} x_k$ is $12$, we have
$\sum_{k \geq 3} x_k (5-k) = 0$. If $P$ has triangular or
quadrangular facets, then $\beta^{-2,8}(P) \neq 0$ by
Theorem~\ref{theorem: Hochster}. Therefore, $x_3 = x_4 = 0$.
Now if $x_k \neq 0$ for $k \geq 6$,
then $\sum_{k \geq 3} x_k (5-k)$ must be negative. This implies
$x_5 = 12$. Hence $P$ is a dodecahedron.
\end{proof}

\section{Some variations of the definition of rigidity}
\enlargethispage{2\baselineskip} There are several variations of
the definition of cohomological rigidity. As is mentioned in the
Introduction, cohomological rigidity is first introduced in
\cite{MS} in terms of toric manifolds and simplical complexes.
Namely, a simplicial complex $\Sigma_X$ associated with a toric
manifold $X$ is {\em rigid} if $\Sigma_X\approx\Sigma_Y$ whenever
$H^\ast(X)\cong H^\ast(Y)$ as graded rings. Therefore our
definition is a variation of the original definition of rigidity.

Moreover we may consider cohomological rigidity of simple convex
polytopes in terms of small covers,
which gives another variation of the definition. Namely, we may
replace \lq quasitoric manifolds' by \lq small covers' and \lq
integral cohomology rings' by \lq mod 2 cohomology rings' in
Definition~\ref{definiton:rigidity}.
A small cover is a closed $n$-dimensional manifold with a locally
standard mod~$2$ torus $(\mathbb Z_2)^n$-action over a simple
convex polytope. It is therefore a mod $2$ analogue of a
quasitoric manifold.
%
Small covers are introduced by Davis and Januszkiewicz
in~\cite{DJ}.

In the proof of our rigidity results we made essential use of bigraded Betti numbers which are
purely combinatorial invariants of the polytopes. Considering this, Buchstaber asked the following
question in his lecture
notes~\cite{Bu}.
\begin{question}
Let $K$ and $K'$ be simplicial complexes, and let $\mathcal Z_K$ and $\mathcal Z_{K'}$ be
their respective moment angle complexes.
When does a cohomology ring isomorphism $H^\ast(\mathcal Z_K : k)\cong H^\ast(\mathcal Z_{K'} : k)$
imply a combinatorial equivalence $K\approx K'$ where $k$ is field?
\end{question}
Let us call the simplicial complexes giving the positive answer to
the question {\em B-rigid}. Note that $H^\ast(\mathcal Z_K:k)\cong
\Tor(k(K), k)$, see \cite{BP}. Let $K=(\partial P)^\ast$ (resp.
$K'=(\partial P')^\ast$) be the dual of the boundary of a simple
convex polytope $P$ (resp. $P'$).  Let $M$ (resp. $M'$) be a
quasitoric manifold over $P$ (resp. $P'$) such that
$H^\ast(M)\cong H^\ast(M)$. Then by Lemma~\ref{lemma:
Tor(P)=Tor(p')} and the ring isomorphism $H^\ast(\mathcal
Z_K:k)\cong \Tor(k(K), k)$, we have the isomorphism
$H^\ast(\mathcal Z_K :k)\cong H^\ast(\mathcal Z_{K'}:k)$. Hence if
$P$ is cohomologically rigid, then $K$ is B-rigid. Furthermore
Example~\ref{vcex} still gives non B-rigid simplicial complexes.
However at this moment we do not know whether cohomological
rigidity
is equivalent to B-rigidity for simple convex polytopes supporting
quasitoric manifolds.

\enlargethispage{5\baselineskip}


\begin{thebibliography}{99}


\bibitem{BB90}
G. Blind \and R. Blind, `Convex polytopes without triangular
faces', {\em Israel J. Math.}, 71(2): 129--134, 1990.
%
\bibitem{BB92}
G. Blind \and R. Blind, `Triangle-free polytopes with few facets',
{\em Arch. Math. (Basel)}, 58(6): 599--605, 1992.
%
\bibitem{BG}
W. Bruns and  J. Gubeladze, `Combinatorial invariance of
{S}tanley-{R}eisner rings', {\em Georgian Math. J.}, 3(4):
315-318, 1996.
%
\bibitem{Bu}
V.\,M. Buchstaber, `Lectures on Toric Topology. lecture notes of
\lq Toric Topology Workshop: KAIST 2008', {\em  Trends in
Mathematics} 10(1) 1--64, 2008.

\bibitem{BP}
V.\,M. Buchstaber \and T.\,E. Panov, `Torus Actions and Their
Applications in Topology and Combinatorics', University Lecture
Series, Vol. 24, Amer. Math. Soc., Providence, R.I., 2002.

\bibitem{CMS}
S. Choi, M. Masuda, \and D.\,Y. Suh, `Topological classification
of generalized Bott towers', {\em preprint}; {\em to appear in
Trans. Amer. Math. Soc.}; {\em arXiv:0807.4334[math.AT]}.


\bibitem{DJ}
M.\,W. Davis \and T. Januszkiewicz, `Convex polytopes, {C}oxeter
orbifolds and torus actions', {\em Duke Math. J.}, 62(2):417--451,
1991.

\bibitem{Ho}
M. Hochster, `Cohen-Macaulay rings, combinatorics, and simplicial
complexes', in: Ring Theory II (Proc. Second Oklahoma Conference),
B. R. McDonald and R. Morris, eds., Dekker, New York, 1977, pp.
171-223.

\bibitem{MP}
M. Masuda \and T. Panov, `Semifree circle actions, Bott towers,
and quasitoric manifolds', {\em Sbornik Mathematics},
199(8):1201--1223, 2008; {\em arXiv:math.AT/0607094v2}.

\bibitem{MS}
M. Masuda \and D. Y. Suh, `Classification problems of toric
manifolds via topology', in: Toric Topology, M.~Harada,
Y.~Karshon, M.~Masuda, and T.~Panov, eds., Contemporary
Mathematics, Vol.~460, Amer. Math. Soc., Providence, R.I., 2007,
pp.~273--286; {\em arXiv:math.AT/0709.4579}.

\bibitem{Od}
T. Oda, `Convex Bodies and Algebraic Geometry. An Introduction to
the Theory of Toric Varieties', Ergeb. Math. Grenzgeb. (3), 15,
Spinger-Verlag, Berlin, 1988.

\end{thebibliography}
\end{document}